\documentclass[12pt]{amsart}

\usepackage{amsmath, amsthm, amssymb}
\usepackage{graphicx}
\usepackage{mathbbol}
\usepackage{txfonts}
\usepackage[usenames]{color}
\usepackage{pgfplots}
\usepackage{comment}

\allowdisplaybreaks[4]

\newtheorem{thm}{Theorem}[section]
\newcommand{\bt}{\begin{thm}}
\newcommand{\et}{\end{thm}}

\newtheorem{cor}[thm]{Corollary}   
\newcommand{\bc}{\begin{cor}}
\newcommand{\ec}{\end{cor}}

\newtheorem{lem}[thm]{Lemma}   
\newcommand{\bl}{\begin{lem}}
\newcommand{\el}{\end{lem}}

\newtheorem{prop}[thm]{Proposition}
\newcommand{\bp}{\begin{prop}}
\newcommand{\ep}{\end{prop}}

\newtheorem{defn}[thm]{Definition}
\newcommand{\bd}{\begin{defn}}    
\newcommand{\ed}{\end{defn}}

\newtheorem{rmrk}[thm]{Remark}   

\newcommand{\br}{\begin{rmrk}}
\newcommand{\er}{\end{rmrk}}

\newcommand{\be}{\begin{equation}}

 \newcommand{\ee}{\end{equation}}

\begin{document}

\title[OBERWOLFACH WORKSHOP REPORT ]
{OBERWOLFACH WORKSHOP REPORT\\
Analysis, Geometry and Topology of \\Positive Scalar Curvature Metrics
\\.\\
{\em Limits of sequences of manifolds with \\
nonnegative scalar curvature and other hypotheses} }

\author{Christina Sormani with Wenchuan Tian and Changliang Wang}
\thanks{Prof. Sormani is partially supported by NSF DMS \#1006059.}
\address{Christina Sormani, CUNY Graduate Center and Lehman College, NY, NY, 10016}
\email{sormanic@gmail.com}


\maketitle



In 2014, Gromov suggested that one should formulate and prove
a compactness theorem for sequences of compact manifolds,
$M_j^m$, with $Scal\ge 0$ and develop a notion of generalized 
$Scal\ge 0$ on the class of limit spaces \cite{Gromov-Plateau}.   He suggested that perhaps
Sormani-Wenger's Intrinsic Flat (SWIF) convergence of \cite{SW-JDG} might work well.
In fact Wenger had already proven a compactness theorem:
$$
Vol(M^m_j)\le V_0 \textrm{ and } Diam(M^m_j)\le D_0 \implies \,\,\exists M_{j_k} \to M_\infty \textrm{ in SWIF sense,}
$$
where $M_\infty$ is an integral current space, possibly the $0$ space \cite{Wenger-compactness}.   
Note that collapsing sequences of round spheres SWIF converge to the $0$ space, so a compactness theorem 
for manifolds with $Scal \ge 0$ would need a noncollapsing hypothesis to avoid such a limit.

In 2018, at an IAS Emerging Topics meeting, a conjecture was formulated in dimension 3 using the two hypotheses
of Wenger's Compactness Theorem and a very weak noncollapsing hypothesis, $MinA(M_j)\ge A_0>0$,
where
$$
MinA(M)=min\{Area(\Sigma)\,: \textrm{ closed min surface }\Sigma\subset M^3\}.
$$
Note that this hypothesis immediately rules out collapsing spheres (whose equators are unstable minimal surfaces)
as well as bubbles (which have necks with stable minimal surfaces).   This was a natural assume
in light of Schoen-Yau theorems regarding stable minimal surfaces, Marques-Neves Theorems regarding
unstable minimal surfaces, and the Penrose Inequality.   See the 2023 survey by Sormani \cite{Sormani-survey-23}.

The $MinA$ hypothesis rules out sequences of $M_j^3$ with $Scal\ge 0$ which are built using increasingly thin Schoen-Yau or Gromov-Lawson tunnels  \cite{Schoen-Yau-structure}\cite{Gromov-Lawson-Spin1}.   This includes the sewing examples of Basilio-Dodziuk-Sormani \cite{BDS-sewing} which GH and SWIF converge to a pulled string space 
and examples of Basilio-Sormani \cite{BS-seq} which GH and SWIF converge to a space where a set, $K\in {\mathbb S}^3$, has been identified to a single point.   The key idea in the latter paper is that one can take $M^3_j$ to be $ {\mathbb S}^3$ with an increasingly dense collection of balls located in $K\subset {\mathbb S}^3$ that are replaced by a collection
of increasingly small tunnels that run between each pair of balls.   

In the Basilio-Sormani paper these tunnels were taken to be increasingly short and it is proven via a scruching lemma that the region, $K$, contracts to a single point
or if $K={\mathbb S}^3$ the sequence SWIF converges to the $0$ space.   Dodziuk proved the tunnels
could be taken arbitrarily thin and approaching any length.  More recently Sweeney has proven tunnels can be built with various bounds on Scalar curvature including scalar curvature arbitrarily close to that of a sphere.   Basilio-Kazaras-Sormani took the tunnels developed by Dodziuk in \cite{Dodziuk-tunnels} running between a ball $B_p(\epsilon_j)$ and a ball $B_q(\epsilon_j)$ to have length close to the Euclidean distance between them as points $p,q\subset {\mathbb S}^3\subset {\mathbb E}^4$ and proved this sequence of $M_j$
SWIF converges to $({\mathbb S}^3, d_{{\mathbb E}^4})$ which is a metric space with no midpoints and no geodesics \cite{BKS-no-geod}.   
The $MinA$ hypothesis prevents all these bad examples involving tunneling.

{\bf Open Question 1:} Suppose a sequence, $(M_j^3, g_j)$, with $Scal\ge 0$ converge in some sense to a smooth Riemannian manifold, $(M_0^3, g_0)$.   Does
$(M_0^3, g_0)$ have $Scal\ge 0$?   Gromov and Bamler proved this is true for $C^0$ convergence of the metric tensors when $M^3_j$ are diffeomorphic to one another and to the limit space \cite{Bamler-C} \cite{Gromov-Dirac}.   
SWIF convergence allows one
to study sequences which are not diffeomorphic.  Does this work for SWIF convergence?   I believe the answer should be no if we allow tunneling.  If we take $(M^3_0,g_0)=({\mathbb S}^3, d_{g})$ where $g$ is any Riemannian metric on the sphere such that $g\le g_{\mathbb S}^3$, one can construct an example of a sequence of $(M_j^3,g_j)$ with $Scal\ge 0$
using tunneling similar to the Basilio-Kazaras-Sormani sequence in \cite{BKS-no-geod} with tunnels running between a ball $B_p(\epsilon_j)$ and a ball $B_q(\epsilon_j)$ that have length close the Riemannian distance $d_{g_0}$.   It may be possible to prove this sequence converges in the SWIF sense to $(M^3_0,g_0)$ which does not have $Scal\ge 0$.   

{\bf Open Question 2:} If $(M_j^3,g_j)$ have $Scalar \ge 0$ and satisfy the $MinA$ hypothesis and converge in the SWIF sense to smooth $(M_0,g_0)$, does $(M_0,g_0)$ have $Scal\ge 0$?   The $MinA$ hypothesis prevents the use of tunneling to construct counterexamples.   If a counter example is found, then one needs to consider a stronger hypothesis than the $MinA$ in the compactness conjecture.   Proving this open question in the affirmative would be very challenging.  An easier case to check would be to assume that the sequence converges in the VADB sense defined by Allen-Perales-Sormani (where it was proven that VADB implies SWIF convergence) \cite{VADB}.

There has been some interesting progress on the IAS $MinA$ Compactness Conjecture described above
by Park-Tian-Wang, Tian-Wang, and Kazaras-Xu in \cite{Park-Tian-Wang-18} \cite{Tian-Wang-compactness} \cite{Kazaras-Xu-drawstrings}.   
All three teams consider fixed manifolds, $M^3$, with
varying metric tensors, $g_j$, of the following forms respectively:
$$
({\mathbb S}^3, dt^2+f_j(t)^2 g_{{\mathbb S}^2}), \,\, 
({\mathbb S}^2\times {\mathbb S}^1,g_{{\mathbb S}^2}+f_j(u)^2 d\theta^2), \textrm{ and } 
({\mathbb T}^2\times {\mathbb S}^1,h_j+f_j(u)^2 d\theta^2).
$$
All three achieve $W^{1,p}$ convergence for $p<2$ to a limit metric tensor, $g_\infty$, of the same form satisfying some distributional notion of $Scal\ge 0$.   Park-Tian-Wang 
prove their $f_\infty$ are continuous and take values in $[0,\infty)$.  Although $S=f_\infty^{-1}(0)$ may be an open set,
they prove SWIF convergence of their $(M,g_j)$ to the metric completion of $(M\setminus S, g_\infty)$ because the disappearing regions are wells.

{\bf Open Question 3}: Recall that we already know there is a SWIF limit space, $(M^3, d_\infty)$, which is a rectifiable metric space (possibly the zero space) in the setting of the IAS $MinA$ Compactness Conjecture.   If $(M^3, g_j) \to (M^3, g_\infty)$ have $g_j \to g_\infty$ in the $W^{1,p}$ sense can one prove 
the SWIF limit $(M^3, d_{\infty})$ is isometric to the metric completion of $(M^3\setminus S, g_\infty)$
where $S$ is the singular set where $g_\infty$ is infinite valued or degenerate?  See work
of Allen-Sormani and Allen-Bryden demonstrating how different these can be without
$Scal\ge 0$ and the $MinA$ hypotheses.   With these hypotheses,  
Sormani-Tian-Wang have an extreme limit space in \cite{STW-extreme} with fibres streching to infinite length
and Kazaras-Xu have announced an example in \cite{Kazaras-Xu-drawstrings}
with a pulled thread limit which should be studied closely.   Are there additional hypothesis that can 
guarantee they are isometric? If these spaces are not isometric, how are they related?   Can one use distributional 
$Scal\ge 0$ on the $W^{1,p}$ limit to say something about the SWIF limit?

The extreme limit space found by Sormani-Tian-Wang is a limit space achieved in the Tian-Wang setting 
where $g_j \le g_{j+1}$ converges in the $W^{1,p}$ sense for $p\in [1,2)$ to a limit tensor, $g_\infty$, that is warped by a function,
 $f_\infty$, which is $\infty$ along two circular fibres of infinite length.   They prove this sequence satisfies the hypotheses of the IAS $MinA$ conjecture.   In upcoming work of Sormani-Tian we will prove  the SWIF limit $(M^3, d_{\infty})$ is isometric to the metric completion of $(M^3\setminus S, g_\infty)$ and also isometric to the GH limit.  So this is an interesting space to study possible notions of generalized $Scal\ge 0$.   Open questions in this direction will appear in this upcoming paper
 including possible generalizations of various such notions that might be tested on this extreme space.   Tian-Wang have already generalized the Lee-LeFloch distributional $Scal\ge 0$ to the limits achieved in \cite{Tian-Wang-compactness}
 \cite{Lee-LeFloch}.
 
It is possible that some very natural notions of generalized $Scal\ge 0$ may fail to hold on some of these $W^{1,p}$ and SWIF limit spaces.   It is possible that a stronger hypothesis than the $MinA$ hypothesis is needed in order to obtain better control on the limit space and stronger convergence than simply SWIF convergence.
  
{\bf Open Question 4:} We could consider ${MinL}(M)$, which is the length of shortest closed geodesic in a closed min surface in $M$.   By a theorem of Croke, ${MinL}(M)\ge L_0>0$ is a stronger hypothesis
than ${MinA}(M)\ge A_0>0$ \cite{Croke-area}.   This $MinL$ hypothesis also provides a lower bound on the Gromov filling volumes of the minimal surfaces in $M$ \cite{Gromov-filling}.    Filling volumes are closely related to SWIF distances and can be used to prevent convergence to $0$ and disappearing regions as seen in work of Sormani-Wenger and Portegies-Sormani \cite{SW-CVPDE}\cite{PS-properties}.  One might also consider bounding the systole, $Sys(M)$,  which is the length of shortest closed geodesic in $M$.  See related work of Nabutovsky, Rotman, and Sabourau \cite{NR-GAFA-06}
\cite{NRS-Sweepouts}.

{\bf Open Question 5:} We could consider a noncollapsing condition defined using constant mean curvature surfaces or isoperimetric regions.   Portegies, Jauregui-Lee, and Jauregui-Lee-Perales have theorems about SWIF converging sequences of manfolds whose volumes converge which may be useful to apply in combination with such an hypothesis \cite{Portegies-evalue}\cite{Jauregui-Lee-VF}\cite{JPP-capacity}.   

{\bf Open Question 6:} Is there a notion of convergence for sequences of distinct Riemannian manifolds, $(M_j,g_j)$, which implies this volume preserving SWIF convergence and also captures the information encoded in a $W^{1,p}$ limit?   Perhaps a notion where the convergence on good diffeomorphic regons is in the $W^{1,p}$ sense and the bad regions have volume decreasing to $0$.   Can one define a distance between Riemannian manifolds which captures this notion?   Exploration in this direction involves a deeper understanding of these examples mentioned above and other possible limit spaces with even more badly behaved singular sets.

Complete citations are below but no room was available in the published Oberwolfach Report.

\bibliographystyle{plain}

\bibliography{2023-Sormani.bib}

\end{document}